\documentclass[leqno]{article}
\usepackage[english]{babel}
\usepackage[cp850]{inputenc}
\usepackage{amscd}
\usepackage{amsthm}
\usepackage{amssymb}
\usepackage{amsmath}
\usepackage{verbatim}
\usepackage{color}
\newtheorem{mydef}{Definition}
\newtheorem{myexa}{Example}

\newtheorem{mytheo}{Theorem}
\newtheorem{myque}{Question}
\newtheorem{mycor}{Corollary}

\newtheorem*{myrem}{{\em Remark}}
\newtheorem*{myprem}{{\em Remarks}}

\newcommand{\pt}{\mbox{$\succ$\hspace{-1ex}$\longrightarrow$}}

\begin{document}
\thispagestyle{empty} \vspace*{1cm}
\begin{center}
{\bf\large }
 \end{center}

\begin{center}
{\sc Conditional Expectation of a Markov Kernel Given Another\\ with Some  Applications in Statistical Inference and Disease Diagnosis}\vspace{2ex}\\
A.G. Nogales\vspace{2ex}\\
Dpto. de Matem\'aticas, Universidad de Extremadura\\
Avda. de Elvas, s/n, 06006--Badajoz, SPAIN.\\
e-mail: nogales@unex.es\\
Fax number: +34 924272911
\end{center}
\vspace{.4cm}
\hspace{\parindent}  {\sc Abstract.}

\noindent
Markov kernels play a decisive role in probability  and mathematical statistics theories, and are 
an extension of the concepts of sigma-field and statistic. Concepts such as independence, sufficiency, completeness, ancillarity or conditional distribution  have been extended previously to Markov kernels. 
In this paper, the concept of conditional expectation of a Markov kernel given another is introduced, setting its first properties. 
An application to clinical diagnosis is provided, obtaining {an} optimality property of the predictive values of a diagnosis test. In a statistical framework, this new probabilistic tool is used to extend to Markov kernels the theorems of {Rao-Blackwell} and Lehmann-Scheff\'e. 
{A result about the completeness of a sufficient statistic is obtained in passing by properly enlarging the family of probabilities.}
As a final statistical scholium, a generalization of a result about the completeness of the family of nonrandomized estimators is given.

\vfill
\begin{itemize}
\item[] \hspace*{-1cm} {\em AMS Subject Class.} (2010): {\em Primary\/}
60E05
{\em Secondary\/} 62F10, 62P10, 62B05
\item[] \hspace*{-1cm} {\em Key words and phrases:} Markov kernel, conditional expectation, clinical diagnosis, 
unbiased estimation.
\item[] \hspace*{-1cm} {\em Short running title:} Conditional expectations for Markov kernels

 \end{itemize}

\newpage

\section{Introduction}

Markov kernels (also referred to as stochastic kernels or transition
probabilities) play an important role in probability theory and
mathematical statistics. Indeed, the conditional distribution of one
random variable given another is a Markov kernel (here, we use the
term random variable as being synonymous of a measurable function
between two arbitrary measurable spaces). In fact, as we shall show
below, every Markov kernel is the conditional distribution of some
random variable given another. A transition matrix in Markov chains theory defines a Markov kernel.
Sampling probabilities and posterior distributions in Bayesian
inference are Markov kernels. In statistical decision theory, randomized procedures (also named
decision rules or, even, strategies) are Markov kernels, while
nonrandomized procedures are statistics. It is well known that, in
some situations, the optimum procedure is a randomized one: for
example, the fundamental lemma of Neyman and Pearson shows how
randomization is necessary to obtain a most powerful test; Lehmann
(2005) also describes many other statistical situations where the
use of randomization is properly justified. Pfanzagl (1994, Example
4.2.2) shows a testing problem where there is no nonrandomized test
at least as good as a certain randomized test.

A Markov kernel can also be considered as a generalization of the
concepts of $\sigma$-field and random variable (or statistic, in a
statistical framework). 

Well known concepts of the theory of probabilities or mathematical statistics, such as independence, completeness, ancillarity or conditional distribution  have been
extended to Markov kernels in Nogales (2013a) and Nogales (2013b).
The reader is referred to Heyer (1982) for the corresponding
extension to the concept of sufficiency in the context of informativity for statistical experiments. Notice that an extension to Markov kernels of the concepts and results of probability and mathematical statistics should not be considered useless, as it is not the extension to Markov kernels (or transitions) of the classical theorems of the product measure and Fubini: it is the version for Markov kernel of this theorems what we need to describe the joint distribution of two random variables $X$ and $Y$ in terms of the marginal distribution of $X$ and the conditional distributions of $Y$ given a value of $X$.  

On the other hand, the conditional expectation $E(Y|X)$ of an integrable $n$-dimensional random variable $Y$ given a random variable $X$ is the main tool in the study of the relationship between them; in fact, $y=E(Y|X=x)$ is the so-called regression curve (in a wide sense) of $Y$ on $X$. Basic properties and results on conditional expectations can be found in almost every graduate text in probability theory, after its mathematical introduction in Kolmogorov (1933). 

In this paper, we introduce a new probabilistic tool: the conditional expectation for Markov kernels. Its relationship with the concept of conditional distribution of a Markov kernel given another is established. Some basic properties, two  examples of calculation of such a conditional expectation, and a representation theorem in terms of conditional expectation for random variables are also given. One of the given examples is applied to clinical diagnosis, where some expectations and conditional expectations for Markov kernels get a specific meaning. We obtain in particular {an} optimality property of the predictive values of a diagnosis test as the point that minimizes two naturally weighted distances to the correct decisions on the subpopulations of ill and non-ill individuals. As far as we know, this interpretation of the predictive values appears here for the first time.

As {a} statistical application, in this paper we make use of such tools to extend to Markov kernels the theorems of {Rao-Blackwell} and Lehmann-Scheff\'{e}. These well known theorems are major milestones of mean unbiased estimation theory, going back to Rao (1945) and Blackwell (1947) regard to the Theorem of Rao-Blackwell, and to Lehmann and Scheff\'{e} (1950) regard to the Theorem of Lehmann-Scheff\'{e}. The reader is referred to Pfanzagl (1994, p. 105) for a version for statistics of these theorems; it is assured even there that a more general version of the Rao-Blackwell theorem can be proved in the same way for randomized estimators when a sufficient and complete statistics exists. In this paper both results are generalized for randomized estimators when a sufficient and complete Markov kernel is known.  Two examples of sufficient Markov {kernels} associated to any statistic are given; two more similar examples are provided for complete Markov kernels. 

Notice also that, as the conditional expectation of a Markov kernel given another is a statistic,  a generalization of a result about the completeness of the family of nonrandomized estimators is also obtained.

{Finally, looking for an elusive example of application of this generalized version of the Lehmann-Scheffé Theorem, a result about getting the completenes of a given sufficient statistic has been obtained. The key:   to enlarge the family of probabilities adequately.}

For ease {of} reading, the proofs of the results have been collected in the last section of the paper. 

\section{Basic definitions}
{The concepts presented in this section can be found in Heyer (1982) (see also Dellacherie and Meyer (1988)) and therefore they will be exposed very briefly, even at risk of being somewhat dense.  However the usual  notations in this area has been modified. 
As explained in the Remak 2 below, the concept of  Markov kernel is an extension of the concept of random variable (and also of the concept of $\sigma$-field) and the notation to be used for operations with Markov kernels,  the same that for random variables, tries to highlight this analogy. }

In the next,  $(\Omega,\mathcal A)$, $(\Omega_1,\mathcal A_1)$,
and so on, will denote measurable spaces. A random variable is a map
$X:(\Omega,\mathcal A)\rightarrow (\Omega_1,\mathcal A_1)$ such that
$X^{-1}(A_1)\in\mathcal A$, for all $A_1\in\mathcal A_1$. Its
probability distribution (or, simply, distribution) $P^X$ with
respect to a probability measure $P$ on $\mathcal A$ is the image
measure of $P$ by $X$, i.e., the probability measure on $\mathcal
A_1$ defined by $P^X(A_1):=P(X^{-1}(A_1))$. We will write $\times$
instead of $\otimes$ for the product of $\sigma$-fields or measures. $\mathcal R^k$ will denote the Borel $\sigma$-field on $\mathbb R^k$.

\begin{mydef}\rm  (Markov kernel) A Markov kernel
$M_1:(\Omega,\mathcal A)\pt   (\Omega_1,\mathcal A_1)$ is a map $M_1:\Omega\times\mathcal A_1\rightarrow[0,1]$ such that\\
(i) $\forall \omega\in\Omega$, $M_1(\omega,\cdot)$ is a  probability
measure on
$\mathcal A_1$,\\
(ii) $\forall A_1\in\mathcal A_1$, $M_1(\cdot,A_1)$ is $\mathcal A$-measurable.
\end{mydef}

\begin{myprem}\rm
1) Given two random variables $X_i:(\Omega,\mathcal
A,P)\rightarrow(\Omega_i,\mathcal A_i)$, $i=1,2$, the conditional
distribution of $X_2$ given $X_1$, when it exists, is a Markov
kernel $M:(\Omega_1,\mathcal A_1)\pt (\Omega_2,\mathcal A_2)$ such
that $P(X_1\in A_1,X_2\in
A_2)=\int_{A_1}M(\omega_1,A_2)dP^{X_1}(\omega_1)$, for all
$A_1\in\mathcal A_1$ and $A_2\in\mathcal A_2$. We write
$P^{X_2|X_1=\omega_1}(A_2):=M(\omega_1,A_2)$. Reciprocally, every
Markov kernel is a conditional distribution; namely, given a Markov
kernel $M_1:(\Omega,\mathcal A,P)\pt (\Omega_1,\mathcal A_1)$, it is
easily checked that 
$$M_1(\omega,A_1)=(P\otimes
M_1)^{\pi_1|\pi=\omega}(A_1),$$
where
$\pi:\Omega\times\Omega_1\rightarrow \Omega$ and
$\pi_1:\Omega\times\Omega_1\rightarrow \Omega_1$ are the coordinatewise projections
and $P\otimes M_1$ stands for the only probability measure on the
product space $(\Omega\times\Omega_1,\mathcal A\times\mathcal A_1)$
such that $(P\otimes M_1) (A\times
A_1)=\int_AM_1(\omega,A_1)\,dP(\omega)$ for all $A\in\mathcal A$ and
$A_1\in\mathcal A_1$.

2) The concept of Markov kernel extends the concepts of random variable
 and $\sigma$-field. A random variable $T_1:(\Omega,\mathcal
A)\rightarrow(\Omega_1,\mathcal A_1)$ will be identified with the
Markov kernel $M_{T_1}:(\Omega,\mathcal A)\pt   (\Omega_1,\mathcal
A_1)$ defined by
$$M_{T_1}(\omega,A_1)=\delta_{T_1(\omega)}(A_1)=I_{A_1}(T_1(\omega)),
$$
where $\delta_{T_1(\omega)}$ denotes the Dirac measure -the
degenerate distribution- at the point $T_1(\omega)$ and $I_{A_1}$ is
the indicator function of the event $A_1$. The sub-$\sigma$-field
$\mathcal B\subset\mathcal A$ will be identified with the Markov
kernel $M_{\mathcal B}:(\Omega,\mathcal A)\pt   (\Omega,\mathcal B)$
given by $M_{\mathcal B}(\omega,B)=\delta_{\omega}(B)$.\end{myprem}

\begin{mydef}\rm
(Image of a Markov kernel) The image (or {probability
distribution}) of a Markov kernel $M_1:(\Omega,\mathcal A,P)\pt
(\Omega_1,\mathcal A_1)$ on a probability space is the probability
measure  $P^{M_1}$ on $\mathcal A_1$ defined by
$$P^{M_1}(A_1):=\int_{\Omega}M_1(\omega,A_1)\,dP(\omega).$$
\end{mydef}

\begin{myrem}\rm  Note that
$$P^{M_1}=(P\otimes M_1)^{\pi_1}
$$
where $\pi_1:\Omega\times\Omega_1\rightarrow\Omega_1$ denotes the
coordinatewise projection onto $\Omega_1$. So, if $f:(\Omega_1,\mathcal A_1)\rightarrow \mathbb
R$ is a nonnegative or $P^{M_1}$-integrable function,
\begin{gather*}\begin{split}
\int_{\Omega_1}f(\omega_1)dP^{M_1}(\omega_1)&=\int_\Omega\int_{\Omega_1}f(\omega_1)M_1(\omega,d\omega_1)dP(\omega)\\
&=\int_{\Omega\times\Omega_1}f(\omega_1)d(P\otimes
M_1)(\omega,\omega_1).
\end{split}\end{gather*}

\end{myrem}

\begin{mydef}\rm  {\rm (a)} (Composition of Markov kernels) The composition of two Markov kernels
$M_1:(\Omega_1,\mathcal A_1)\pt (\Omega_2,\mathcal A_2)$ and
$M_2:(\Omega_2,\mathcal A_2)\pt (\Omega_3,\mathcal A_3)$ is defined
as the Markov kernel $$M_2M_1:(\Omega_1,\mathcal A_1)\pt
(\Omega_3,\mathcal A_3)$$ 
given by
$$M_2M_1(\omega_1,A_3)=\int_{\Omega_2}M_2(\omega_2,A_3)M_1(\omega_1,d\omega_2).
$$

{\rm (b)} (Composition of a Markov kernel and a random variable) Let
$X_1:(\Omega,\mathcal A)\rightarrow (\Omega_1,\mathcal A_1)$ be a
random variable and $M_1:(\Omega_1,\mathcal A_1)\pt
(\Omega'_1,\mathcal A'_1)$ a Markov kernel. A new Markov kernel
$M_1X_1:(\Omega,\mathcal A)\pt (\Omega'_1,\mathcal A'_1)$ is defined
by means of
$$M_1X_1(\omega,A'_1):=M_1(X_1(\omega),A'_1).
$$
\end{mydef}

\begin{myrem}\rm
When $M_{X_1}$ is the Markov kernel corresponding to the random
variable  $X_1$, we have that $M_1X_1=M_1M_{X_1}$.
\end{myrem}

\section{Expectation and conditional expectation for Markov kernels}

{In this section, of a character rather technical, we first recall the concept of expectation of a Markov kernel and introduce the definition of conditional expectation for Markov kernels, the last being parallel to the definition of conditional expectation for random variables that can be found, for instance, in Heyer (1982, p. 264). Some examples of calculation and some possible applications to clinical diagnoses willl be given joint with some properties --even a optimality one-- of this new concept.}

Let $(\Omega,\mathcal A,P)$ be a probability space. {Onwards, $\mathbb R^k$ is supposed to be endowed with its Borel $\sigma$-field $\mathcal R^k$. }

\begin{mydef}\rm  (Expectation of a Markov kernel) A Markov kernel $M_1:(\Omega,\mathcal A,P)\pt \mathbb R^k$ is said to be $P$-integrable if the map $\omega\mapsto \int_{\mathbb R^k}xM_1(\omega,dx)$ is $P$-integrable, i.e., if there exists and is finite the integral
	$$\int_\Omega\int_{\mathbb R^k}xM_1(\omega,dx)dP(\omega)
	$$
or, which is the same, if the distribution $(P\otimes M_1)^{\pi_2}$ has finite mean, where $\pi_2:\Omega\times\mathbb R^k\rightarrow\mathbb R^k$ denotes the second coordinatewise projection. In this case, we define the expectation of the Markov kernel $M_1$ as
	$$E_P(M_1):=\int_\Omega\int_{\mathbb R^k}xM_1(\omega,dx)dP(\omega)
	$$
\end{mydef}

\begin{myrem}\rm
{Notice that the previous definition extends standard expectation for a random variable $X:(\Omega,\mathcal A,P)\rightarrow\mathbb R^k$, as $E_P(M_X)=E_P(X)$. }
\end{myrem}

\begin{mydef}\rm
	Let $M_1:(\Omega,\mathcal A,P)\pt \mathbb R^k$ be a $P$-integrable Markov kernel. We define a set function $M_1\cdot P$ on $\mathcal A$ by
	$$(M_1\cdot P)(A):=\int_A\int_{\mathbb R^k}xM_1(\omega,dx)dP(\omega).
	$$
\end{mydef}

Note that $M_1\cdot P\ll P$ and $(M_1\cdot P)^{M_2}\ll P^{M_2}$, when $M_2:(\Omega,\mathcal A,P)\pt (\Omega_2,\mathcal A_2)$ is another Markov kernel.

\begin{mydef}\rm \label{def8} (Conditional expectation of a Markov kernel given another) Let $M_1:(\Omega,\mathcal A,P)\pt \mathbb R^k$ be a $P$-integrable Markov kernel and $M_2:(\Omega,\mathcal A,P)\pt (\Omega_2,\mathcal A_2)$ be a Markov kernel. The conditional expectation $E_P(M_1|M_2)$ is defined by:
	$$E_P(M_1|M_2):=\frac{d(M_1\cdot P)^{M_2}}{dP^{M_2}}
	$$
	i.e., $E_P(M_1|M_2)$ is the (equivalence class of) real measurable function(s) on $(\Omega_2,\mathcal A_2)$ such that, for every $A_2\in\mathcal A_2$,
	\begin{gather*}\begin{split}
			\int_\Omega M_2(\omega,A_2)\int_{\mathbb R^k}xM_1(\omega,dx)dP(\omega)&=\int_{A_2}E_P(M_1|M_2)dP^{M_2}\\
			&=\int_\Omega\int_{A_2}E_P(M_1|M_2)(\omega_2)M_2(\omega,d\omega_2)dP(\omega).
		\end{split}\end{gather*}
	\end{mydef}
	
	The next result yields an integral representation of such a conditional expectation. First, we refer the reader to Nogales (2013b) for the definition and existence of the conditional distribution $P^{M_1|M_2}$ of a Markov kernel  $M_1:(\Omega,\mathcal A,P)\pt (\Omega_1,\mathcal A_1)$ with respect to another Markov kernel $M_2:(\Omega,\mathcal A,P)\pt (\Omega_2,\mathcal A_2)$. Namely, it is defined as a Markov kernel
	$L:(\Omega_2,\mathcal A_2)\pt   (\Omega_1,\mathcal A_1)$ such that,
	for every pair of events $A_1\in\mathcal A_1$ and $A_2\in\mathcal A_2$,
	\begin{gather*}\begin{split}
			\int_\Omega
			M_1(\omega,A_1)M_2(\omega,A_2)dP(\omega)&=\int_{A_2}L(\omega_2,A_1)dP^{M_2}(\omega_2)\\&=
			\int_\Omega\int_{A_2}L(\omega_2,A_1)M_2(\omega,d\omega_2)dP(\omega)
		\end{split}
	\end{gather*}

	\begin{mytheo}\rm \label{teo1} Let $M_1$ and $M_2$ be two Markov kernels as in the previous definition. Then
		$$E_P(M_1|M_2)(\omega_2)=\int_{\mathbb R^k}xP^{M_1|M_2}(\omega_2,dx)
		$$
		(in the sense that the last integral defines a version of the conditional expectation of $M_1$ given $M_2$). More generally, if $f:\mathbb R^k\rightarrow\mathbb R^m$ has nonnegative components or is $P^{M_1}$-integrable function, then
		$$E_P(fM_1|M_2)(\omega_2)=\int_{\mathbb R^k}f(x)P^{M_1|M_2}(\omega_2,dx),$$
		where $fM_1$ is the Markov kernel defined by $fM_1(\omega,C):=M_1(\omega,f^{-1}(C))$, $\omega\in\Omega$, $C\in\mathcal R^m$.
	\end{mytheo}

The following are two examples of calculation. 

\begin{myexa}\rm Given $\theta\in[0,1]$, let $\Omega=\{0,1\}$, $\mathcal A=\mathcal P(\Omega)$ and $P$  the probability measure  on $(\Omega,\mathcal A)$ assigning probability $\theta$ to the point 1 and $1-\theta$ to the point 0. For $i=1,2$, consider the Markov kernel $M_i:(\Omega,\mathcal A)\pt (\Omega,\mathcal A)$ defined by the stochastic matrix
	$$\left(\begin{array}{cc}
	p_i & 1-p_i \\
	q_i & 1-q_i\\
	\end{array}\right),
	$$
	where $0\le p_i,q_i\le 1$. Then
$$P^{M_1}(\{0\})=\int_{\{0,1\}}M_1(\omega,\{0\})dP(\omega)=(1-\theta)p_1+\theta q_1
$$
and $P^{M_1}(\{1\})=(1-\theta)(1-p_1)+\theta (1-q_1)$. Hence,
$$E_P(M_1)=\int_{\{0,1\}}\int_{\mathbb R}xM_1(\omega,dx)dP(\omega)=
P^{M_1}(\{1\}).$$ 

Moreover, if $L:=P^{M_2|M_1}:(\Omega,\mathcal A)\pt (\Omega,\mathcal A)$, according to Nogales (2013b, Prop. 2), given $\omega_1,\omega_2\in\{0,1\}$, 
$$L(\omega_1,\{\omega_2\})=
\frac{\int_{\{0,1\}}M_1(i,\omega_1)M_2(i,\omega_2)dP(i)}{\int_{\{0,1\}}M_1(i,\omega_1)dP(i)}
=\frac{(1-\theta)M_1(0,\omega_1)M_2(0,\omega_2)+\theta M_1(1,\omega_1)M_2(1,\omega_2)}
{(1-\theta)M_1(0,\omega_1)+\theta M_1(1,\omega_1)}
$$
Hence
$$L(0,\{1\})=\frac{(1-\theta)M_1(0,0)M_2(0,1)+\theta M_1(1,0)M_2(1,1)}
{(1-\theta)M_1(0,0)+\theta M_1(1,0)}=\frac{(1-\theta)p_1(1-p_2)+\theta q_1(1-q_2)}
{(1-\theta)p_1+\theta q_1}
$$
and
$$L(1,\{1\})=\frac{(1-\theta)M_1(0,1)M_2(0,1)+\theta M_1(1,1)M_2(1,1)}
{(1-\theta)M_1(0,1)+\theta M_1(1,1)}=\frac{(1-\theta)(1-p_1)(1-p_2)+\theta (1-q_1)(1-q_2)}
{(1-\theta)(1-p_1)+\theta (1-q_1)},
$$
while $L(\omega_1,\{0\})=1-L(\omega_1,\{1\})$, $\omega_1=0,1$. Finally, for $\omega_1\in\{0,1\}$,
$$E_P(M_2|M_1)(\omega_1)=\int_{\{0,1\}}xP^{M_2|M_1}(\omega_1,dx)=L(\omega_1,\{1\}).
$$
\indent {\sc Subexample 1.1:} (Application to clinical diagnosis) Consider a diagnosis tests $T$ for a certain disease $D$. We write $D=1$ ($=0$) for an individual having (not having) the disease as determined by a ``gold standard" diagnostic procedure, and $T=1$ ($=0$) if the diagnostic is positive (negative). There are several terms that are commonly used in this context: $P(D=1)$ is called the prevalence of the disease (on a given population), while $s=P(T=1|D=1)$ is the sensitivity of the test and $e=P(T=0|D=0)$ is its specificity. The stochastic matrix
$$M_1=\left(\begin{array}{lc}
p_1=e & 1-e \\
q_1=1-s & s\\
\end{array}\right),
$$
describes the transition probabilities from the state $i\in\{0,1\}$ (the gold standard test is negative -$i=0$- or positive -$i=1$-) to the state $j\in\{0,1\}$ (the test $T$ is negative -$j=0$- or positive -$j=1$-). This way, $M_1$ becomes a Markov kernel from $\{0,1\}$ to $\{0,1\}$ and its probability distribution $P^{M_1}$ satisfies
$$P^{M_1}(\{1\})=(1-\theta)(1-e)+\theta(1-s)=P(T=1),
$$
the probability that any given individual of the population receive a positive diagnostic. 
If $M_2$ denotes the gold standard diagnostic test, $M_2$ is (identified with) the identity matrix of order 2 (i.e., $p_2=1$ and $q_2=0$). So, analogously, $P^{M_2}(\{1\})=P(D=1)=\theta$. Moreover, if $L:=P^{M_2|M_1}$, according to the example, we have that
$$E_P(M_2|M_1)(1)=L(1,\{1\})=\frac{\theta s}{(1-\theta)(1-e)+\theta s}=P(D=1|T=1)
$$
is the so-called predictive positive value $PPV$ of the diagnosis test $T$ and, in the same way,
$$E_P(M_2|M_1)(\{0\})=L(0,\{1\})=1-\frac{(1-\theta)e}{(1-\theta)e+\theta(1-s)}=1-P(D=0|T=0),
$$
$P(D=0|T=0)$ being the predictive negative value $PNV$ of $T$. Now, if
$$N_1=\left(\begin{array}{cc}
1-e & e \\
1-s & s\\
\end{array}\right),
$$
we have that
$$E_P(N_1)=P^{N_1}(\{1\})=(1-\theta)e+\theta s=P(T=0,D=0)+P(T=1,D=1)
$$
is the ``accuracy" of $T$, i.e., the proportion of true diagnostics of $T$. Moreover
$$E_P(M_2|N_1)(1)=\frac{\theta s}{\theta s+(1-\theta)e}=\frac{P(T=1,D=1)}{P(T=0,D=0)+P(T=1,D=1)}
$$
is the proportion of positive true diagnostics among all true diagnostics of $T$. $\Box$

\end{myexa}
		
	\begin{myexa}\rm
		For $1\le i\le 3$, let $(\Omega_i,\mathcal A_i,\mu_i)$ be a
		$\sigma$-finite measure space such that $(\Omega_i,\mathcal A_i)$ is
		a standard Borel space for $i=2,3$, and $X_i:(\Omega,\mathcal
		A,P)\rightarrow (\Omega_i,\mathcal A_i,\mu_i)$ is a random variable.
		We assume that the joint distribution of $X=(X_1,X_2,X_3)$ admits a
		density $f$ with respect to the product measure
		$\mu_1\times\mu_2\times\mu_3$. We write $f_{ij}$ for the joint
		$\mu_i\times\mu_j$-density of $(X_i,X_j)$ when $1\le i<j\le 3$, and
		$f_i$ for the $\mu_i$-density of $X_i$.  It is shown in Nogales
		(2013b, Example 1) that the conditional distributions
		$M_i=P^{X_i|X_1}:(\Omega_1,\mathcal A_1)\pt (\Omega_i,\mathcal
		A_i)$, $i=2,3$,  and $L:=P_1^{M_2|M_3}$ exist, where $P_1=P^{X_1}$,
		and that
		a density of $L(\omega_3,\cdot)$ with respect to $\mu_2$ is the map
		$$\omega_2\mapsto\int_{\Omega_1}\frac{f_{12}(\omega_1,\omega_2)f_{13}(\omega_1,\omega_3)}{f_1(\omega_1)f_3(\omega_3)}d\mu_1(\omega_1)
		$$
		$L$ is in fact the conditional distribution of a conditional distribution given another conditional distribution!  So, when $(\Omega_2,\mathcal A_2)=(\mathbb R^k,\mathcal R^k)$ and
		$\mu_2$ is the Lebesgue measure, we have that the conditional
		expectation $E_{P_1}(M_2|M_3)$ is the map
		$$\omega_3 \mapsto \int_{\mathbb R^k}x_2\int_{\Omega_1}\frac{f_{12}(\omega_1,x_2)f_{13}(\omega_1,\omega_3)}{f_1(\omega_1)f_3(\omega_3)}d\mu_1(\omega_1)dx_2
		$$
		For instance, let $X=(X_1,X_2,X_3)$ be a trivariate normal random
		variable with null mean and $P_1$ the marginal distribution of
		$X_1$. For $i=2,3$, consider the Markov kernel $M_i=P_1^{X_i|X_1}$,
		the conditional distribution of $X_i$ given $X_1$. It is shown in
		Nogales (2013b) that the conditional distribution $L:=P_1^{M_2|M_3}$
		of the Markov kernel $M_2$ given $M_3$ with respect to $P_1$ satisfy
		that $L(x_3,\cdot)$ is the univariate normal distribution of mean
		$\frac{\sigma_2\rho_{12}\rho_{13}}{\sigma_3}x_3$ and variance
		$\sigma_2^2(1-\rho_{12}^2\rho_{13}^2)$, where $\sigma_i$ is the
		standard deviation of $X_i$ and $\rho_{ij}$ stand for the
		correlation coefficient of $X_i$ and $X_j$. According to the
		previous result, the conditional expectation of $M_2$ given $M_3$ is
		the random variable $x_3\mapsto
		\frac{\sigma_2\rho_{12}\rho_{13}}{\sigma_3}x_3$.\par 
		Moreover, it is easily checked that $P_1^{M_2}=P^{X_2}$ and $E_{P_1}(M_2)=E_P(X_2)$. $\Box$ \end{myexa}

	Note that
	\begin{gather*}
		\begin{split}
			\int_{\Omega_2}E_P(M_1|M_2)dP^{M_2}&=\int_{\Omega_2}\int_{\mathbb R^k}xP^{M_1|M_2}(\omega_2,dx)dP^{M_2}(\omega_2)\\
			&=\int_{\Omega_2\times\mathbb R^k}xdP^{M_2\times M_1}(\omega_2,x)\\&=
			\int_{\mathbb R^k}xd\left(P^{M_2\times M_1}\right)^{\pi}(x),
		\end{split}\end{gather*}
		where $\pi:\Omega_2\times\mathbb R^k\rightarrow\mathbb R^k$ is the coordinatewise projection and $M_2\times  M_1:(\Omega,\mathcal A)\pt (\Omega_2\times \mathbb R^k,\mathcal A_2\times\mathcal R^k)$ satisfies
		$(M_2\times  M_1)(\omega,A_2\times A_1)= M_2(\omega,A_2)\cdot M_1(\omega,A_1)$, $A_i\in\mathcal A_i$,
		$i=1,2.$ But $\left(P^{M_2\times M_1}\right)^{\pi}=P^{M_1}$. Hence
		$$E_{P^{M_2}}(E_P(M_1|M_2))=\int_{\Omega_2}E_P(M_1|M_2)dP^{M_2}=\int_{\mathbb R^k}xdP^{M_1}(x)=E_P(M_1).
		$$
		
		This way we obtain the following corollary, which generalizes a known property of usual conditional expectations.
		\begin{mycor}\rm \label{cor1} Let $M_1$ and $M_2$ be two Markov kernels as in the previous Definition \ref{def8}. Then
			$$E_{P^{M_2}}(E_P(M_1|M_2))=E_P(M_1).
			$$
		\end{mycor}

		We can have a representation of conditional expectations for Markov kernels in terms of conditional expectations for random variables.
		
		\begin{mytheo}\rm \label{teo2} If $M_1$ is $P$-integrable,
		$E_P(M_1|M_2)=E_{P\otimes M_2}(\bar{M}_1|\pi_2)$ where $\bar
			M_1:(\Omega\times\Omega_2,\mathcal A\times\mathcal
			A_2)\rightarrow\mathbb R^k$ is defined by $\bar
			M_1(\omega,\omega_2):=\int_{\mathbb R^k}xM_1(\omega,dx)$, and
			$\pi_2$ is the second coordinatewise projection on
			$\Omega\times\Omega_2$.
		\end{mytheo}

As a consequence of this representation theorem and Jensen's Inequality, we have the next result. 

\begin{mycor}\rm For every $Z\in\mathcal L^2(\Omega_2,\mathcal A_2,(P\otimes M_2)^{\pi_2}),$ we have that
\begin{gather*}
\|\bar M_1-E_P(M_1|M_2)\|_2^2\le \|\bar M_1-Z\|_2^2,
\end{gather*}
i.e.,
\begin{gather*}\begin{split}
&\int_{\Omega\times\Omega_2}(\bar M_1(\omega,\omega_2)-E_P(M_1|M_2)(\omega_2))^2d(P\otimes M_2)(\omega,\omega_2)\le \\
&\int_{\Omega\times\Omega_2}(\bar M_1(\omega,\omega_2)-Z(\omega_2))^2d(P\otimes M_2)(\omega,\omega_2),
\end{split}\end{gather*}
\end{mycor}

\indent {\sc Subexample 1.1 (cont.):} (Application to clinical diagnosis) Applying the preceding Corollary to Subexample 1.1, writing $a=Z(0)$ ($a$ could represent the probability that the decision 0 is taken, i.e., the test $T$ discards the disease) and $b=Z(1)$ ($b$ could represent the probability that the decision 1 is taken, i.e., the test $T$ confirms the disease), we obtain the following interpretation of predictive values of a diagnostic test $T$:
$$(1-PNV,PPV)=\mbox{arg\,min}_{(a,b)\in\mathbb R^2}\{[(1-a)^2e+b^2(1-e)](1-\theta)+[a^2(1-s)+(1-b)^2s]\theta\}. 
$$
Notice that, for a non-ill individual (i.e., when $D=0$), the right decision will be $(a_0,b_0)=(1,0)$, and $(1-a)^2e+b^2(1-e)$ is a weighted squared distance between $(a,b)$  and the optimal point $(1,0)$ on $\{D=0\}$; the weights are $e=P(T=0|D=0)$ and $1-e=P(T=1|D=0)$ for the discrepancy between $a$ and $a_0=1$,  and that of $b$ and $b_0=0$, respectively, as can be expected. Analogously, for an ill individual (i.e., when $D=1$), the right decision is $(a_1,b_1)=(0,1)$, and $a^2(1-s)+(1-b)^2s$ is also a properly weighted squared distance between $(a,b)$ and the optimal point $(0,1)$ on $\{D=1\}$.

 Notice finally that, in the daily clinical practice, it is not known whether $D=0$ or $D=1$ and we should choose $(a,b)$ in such a way that its simultaneous distance to $(1,0)$ on $\{D=0\}$ and to $(0,1)$ on $\{D=1\}$ reach a minimum; obviously, this simultaneous squared distance is weighted according to the sizes of the subpopulations $\{D=0\}$ and $\{D=1\}$. $\Box$

\section{Some statistical applications: extension to Markov kernels of the Rao-Blackwell and the Lehmann-Scheff\'{e} theorems}

Now, we position ourselves in a statistical context. Let $(\Omega,\mathcal A,\mathcal P)$ be a statistical experiment (i.e., $\mathcal P$ is a family of probability measures on the measurable space $(\Omega,\mathcal A)$). 

The theorems of Rao-Blackwell and Lehmann-Scheff\'e are central results of unbiased point estimation theory. We pursue in this section a version in the Markov kernel framework. 

The concepts defined in the preceding sections can be extended to a statistical framework in a standard way. The concept of sufficiency for Markov kernels is introduced in Heyer (1982, p.163). Recall that, given a Markov kernel $M_1:(\Omega,\mathcal A,\mathcal P)\pt
(\Omega_1,\mathcal A_1)$ and $P\in\mathcal P$, the conditional
probability $P(A|M_1)$ of an event $A\in\mathcal A$ given $M_1$ is
defined as the Radon-Nikodym derivative $d(I_A\cdot
P)^{M_1}/dP^{M_1}$, where $I_A\cdot P$ denotes the measure defined
on $\mathcal A$ by $(I_A\cdot P)(B)=P(A\cap B)$. In other words,
$P(A|M_1)$ is the (equivalence class of) real random variable(s) on
$(\Omega_1,\mathcal A_1)$ such that, for every $A_1\in\mathcal A_1$,

\begin{gather}\begin{split}
\int_A M_1(\omega,A_1) dP(\omega)&=\int_{A_1}
P(A|M_1)dP^{M_1}\\&=\int_\Omega\int_{A_1}P(A|M_1)(\omega_1)
M_1(\omega,d\omega_1)dP(\omega)
\end{split}\end{gather}

\begin{mydef}\rm
(Sufficiency of a Markov kernel) A Markov kernel
$M_1:(\Omega,\mathcal A,\mathcal P)\pt   (\Omega_1,\mathcal A_1)$ is
said to be sufficient if, for  every $A\in\mathcal A$, there exists a
common version $f_A:(\Omega_1,\mathcal A_1)\rightarrow [0,1]$ to the
conditional probabilities $P(A|M_1)$, $P\in\mathcal P$.\end{mydef}

\begin{myprem}\rm
1) The previous definition generalizes that of a sufficient
statistic in the sense that a statistic $T_1$ is sufficient {if and
only if} the corresponding kernel $M_{T_1}(\omega,
A_1)=\delta_{T_1(\omega)}(A_1)$ is sufficient. Also, a
sub-$\sigma$-field $\mathcal B\subset\mathcal A$ is sufficient {if and
	only if} its corresponding kernel $M_{\mathcal B}:(\Omega,\mathcal
A)\pt   (\Omega,\mathcal B)$, defined by $M_{\mathcal
B}(\omega,B):=\delta_\omega(B)$, is also.\par

2) Theorem 22.3 of Heyer (1982) shows that a Markov kernel
$M_1:(\Omega,\mathcal A,\mathcal P)\pt  (\Omega_1,\mathcal A_1)$ is
sufficient { if and	only if} the $\sigma$-field $\pi_1^{-1}(\mathcal
A_1)$ is sufficient in the statistical experiment
$(\Omega\times\Omega_1,\mathcal A\times\mathcal A_1,\{P\otimes
M_1\colon P\in\mathcal P\})$, where $\pi_1$ denotes the coordinatewise projection
over $\Omega_1$.

3) (Sufficiency of Markov kernels when densities are available) Suppose that $\mathcal P$ is dominated by a $\sigma$-finite measure $\mu$ on $(\Omega,\mathcal A)$ --$\mu$ is { typically} the Lebesgue measure in the absolute continuous case and the counting measure in the discrete case--. Let $f_P$ be a $\mu$-density of $P\in\mathcal P$. Let $M_1:(\Omega,\mathcal A,\mathcal P)\pt  (\Omega_1,\mathcal A_1)$ be a Markov kernel and suppose that $m_1:(\Omega\times\Omega_1,\mathcal A\times\mathcal A_1)\rightarrow[0,\infty[$ is a measurable function such that, for every $\omega\in\Omega$, $m_1(\omega,\cdot)$ is a $\mu_1$-density of the probability measure $M_1(\omega,\cdot)$, where $\mu_1$ is a $\sigma$-finite measure on $(\Omega_1,\mathcal A_1)$. It is readily shown that
$$\frac{d(P\otimes M_1)}{d(\mu\times\mu_1)}(\omega,\omega_1)=m_1(\omega,\omega_1)\cdot f_P(\omega).
$$
According to the previous remark and the factorization theorem, the Markov kernel $M_1$ is sufficient {if and	only if} there exist a measurable function $h:(\Omega\times\Omega_1,\mathcal A\times\mathcal A_1)\rightarrow[0,\infty[$ and, for each $P\in\mathcal P$, a measurable function $g_P:(\Omega_1,\mathcal A_1) \rightarrow[0,\infty[$ such that
$$m_1(\omega,\omega_1)\cdot f_P(\omega)=g_P(\omega_1)\cdot h(\omega,\omega_1),\quad\forall \omega,\omega_1.
$$
\end{myprem}

Here we introduce two examples, one discrete and one continuous, of sufficient Markov kernels not associated to statistics. 

\begin{myexa}\rm  Let $\Omega=\{1,2,3\}$, $\mathcal A$ the $\sigma$-field of all subsets of $\Omega$, and $\mathcal P:=\{P_\theta\colon \theta\in[0,1]\}$, where $P_\theta$ assigns probability $\theta/3$ to the points 1 and 2 and probability $1-2\theta/3$ to the point 3. The Markov kernel $M:(\Omega,\mathcal A)\pt (\Omega,\mathcal A)$ defined by the stochastic matrix
$$\left(\begin{array}{ccc}
1/3 & 2/3 & 0\\
1/3 & 2/3 & 0\\
0 & 0 &1
\end{array}\right)
$$
is sufficient and is not associated to any statistic. $\Box$
\end{myexa}

\begin{myexa}\rm\label{continuousexample}  Let $(\Omega,\mathcal A)=(\mathbb R^+,\mathcal R^+)$ and $\mathcal P=\{P_\theta\colon \theta=0,1,2,\dots\}$, where $dP_\theta(x)=I_{[\theta,\theta+1[}(x)\,dx$. For $x\ge 0$, we denote by $M(x,\cdot)$ the uniform distribution on the interval $[\lfloor x\rfloor,\lfloor x\rfloor+1[$, where $\lfloor x\rfloor$ stands for the integer part of $x$. The Markov kernel $M:(\Omega,\mathcal A)\pt (\Omega,\mathcal A)$ is sufficient and is not associated to any statistic. 
	$\Box$
\end{myexa}

Let us recall from Nogales (2013a) the generalization of the concept of completeness to
Markov kernels.

\begin{mydef}\rm  (Completeness of Markov kernels) \label{def8} A Markov kernel $M_1:(\Omega,\mathcal A,\mathcal P)\pt
(\Omega_1,\mathcal A_1)$ is said  to be complete (respectively,
boundedly complete) if, for every (respectively, bounded) real
statistic $f:(\Omega_1,\mathcal A_1,\{P^{M_1}\colon P\in\mathcal
P\})\rightarrow\mathbb R$,
$$E_{P^{M_1}}f=0,\;\forall P\in\mathcal P\quad\Longrightarrow\quad f=0,\;
P^{M_1}\mbox{-almost surely }, \;\forall P\in\mathcal P.
$$
\end{mydef}

\begin{myprem}\rm  1) A Markov kernel
$M_1:(\Omega,\mathcal A,\mathcal P)\pt   (\Omega_1,\mathcal A_1)$ is (respectively, boundedly) complete {if and	only if}  the
$\sigma$-field $\pi_1^{-1}(\mathcal A_1)$ on the statistical
experiment $(\Omega\times\Omega_1,\mathcal A\times\mathcal
A_1,\{P\otimes M_1\colon P\in\mathcal P\})$ is also, where $\pi_1$
denotes the coordinatewise projection over $\Omega_1$, which in turn is equivalent
to the (bounded) completeness of $\pi_1$ (see Nogales  (2103a)). Moreover, if $M_1$ is the Markov kernel corresponding to a statistic $T_1$, then $M_1$ is (boundedly) complete {if and
	only if} $T_1$ is also.

2) (Completeness of Markov kernels when densities are available) Suppose that $\mathcal P$ is dominated by a $\sigma$-finite measure $\mu$ on $(\Omega,\mathcal A)$. Let $f_P$ be a $\mu$-density of $P\in\mathcal P$. Let $M_1:(\Omega,\mathcal A,\mathcal P)\pt  (\Omega_1,\mathcal A_1)$ be a Markov kernel and suppose that $m_1:(\Omega\times\Omega_1,\mathcal A\times\mathcal A_1)\rightarrow[0,\infty[$ is a measurable function such that, for every $\omega\in\Omega$, $m_1(\omega,\cdot)$ is a $\mu_1$-density of the probability measure $M_1(\omega,\cdot)$, where $\mu_1$ is a $\sigma$-finite measure on $(\Omega_1,\mathcal A_1)$. It is readily shown that
$$\frac{d(P\otimes M_1)}{d(\mu\times\mu_1)}(\omega,\omega_1)=m_1(\omega,\omega_1)\cdot f_P(\omega).
$$
According to the previous remark, the Markov kernel $M_1$ is complete if and only if for every statistic $f:(\Omega_1,\mathcal A_1)\rightarrow\mathbb R$ we have that
$$\int_{\Omega\times\Omega_1}f(\omega_1)m_1(\omega,\omega_1) f_P(\omega)d(\mu\times\mu_1)(\omega,\omega_1)=0,\  \ \forall P\in\mathcal P\quad\Longrightarrow\quad
f=0,\ \ (P\otimes M_1)^{\pi_1}-c.s.,\ \forall P\in\mathcal P.
$$
\end{myprem}

Here we present two examples of complete Markov kernels not associated to statistics. 

\begin{myexa}\rm  Let $\Omega=\{1,2\}$, $\mathcal A=\mathcal P(\Omega)$ and $\mathcal P:=\{P_\theta\colon \theta\in[0,1]\}$, where $P_\theta$ assigns probability $\theta$ to the point 1 and probability $1-\theta$ to the point 2. The Markov kernel $M:(\Omega,\mathcal A)\pt (\Omega,\mathcal A)$ defined by the stochastic matrix
$$\left(\begin{array}{cc}
p & 1-p \\
q & 1-q\\
\end{array}\right)
$$
is complete for $p,q\in[0,1]$ when $p\ne q$, and it is not associated to any statistic unless $p,q\in\{0,1\}$. $\Box$
\end{myexa}

\begin{myexa}\rm  Let $\Omega=\mathbb R^+$, $\mathcal A=\mathcal R^+$ and $\mathcal P:=\{P_\theta\colon \theta>0\}$, where $P_\theta$ denotes the exponential distribution of parameter $\theta$. For $x>0$, we denote by $M(x,\cdot)$ the uniform distribution on the interval $[x,x+1[$. The Markov kernel $M:(\Omega,\mathcal A)\pt (\Omega,\mathcal A)$ is complete and is not associated to any statistic. $\Box$
\end{myexa}

Now we are ready to obtain a first extension to Markov kernels of the
theorem of Lehmann-Scheff\'{e}. Theorem \ref{teo6} yields a more general result. First, recall that an statistic $T:(\Omega,\mathcal A,\mathcal P)\rightarrow \mathbb R^k$ is said to be an
unbiased estimator of a function $f:\mathcal P\rightarrow \mathbb R^k$ whenever $E_P(T)=f(P)$, for all $P\in\mathcal P$. $T$ is said to be a minimum variance estimator of $f$ if it is unbiased and has less variance than any other unbiased estimator of $f$.

{ Let $M_1:(\Omega,\mathcal A,\mathcal P)\pt (\Omega_1,\mathcal A_1)$ a Markov kernel and $T:(\Omega,\mathcal A,\mathcal P)\rightarrow \mathbb R^k$ be a statistic. We say that $T$ is a measurable function of $M_1$ if  there exists
a measurable map $S:(\Omega_1,\mathcal A_1)\rightarrow\mathbb R^k$ such that $M_T=M_SM_1$.}

\begin{mytheo}\rm \label{teo3}  Assuming the previous notations, let us suppose that  $M_1:(\Omega,\mathcal A,\mathcal P)\pt (\Omega_1,\mathcal A_1)$ is a sufficient and complete Markov kernel and $T:(\Omega,\mathcal A,\mathcal P)\rightarrow \mathbb R^k$ is an unbiased estimator of a function $f:\mathcal P\rightarrow \mathbb R^k$. If $T$ is a measurable function of $M_1$, then it is the minimum variance unbiased estimator of $f$.
\end{mytheo}

Now let us recall  the definition of unbiased (randomized) estimator.

\begin{mydef}\rm  {\rm (Unbiased estimator)}
An unbiased estimator of a function $f:\mathcal P\rightarrow \mathbb
R^k$ is a $\mathcal P$-integrable Markov kernel $M:(\Omega,\mathcal A,\mathcal P)\pt (\mathbb R^k,\mathcal R^k)$
such that
$$E_P(M):=\int_{\Omega}\int_{\mathbb
R^k}xM(\omega,dx)dP(\omega)= f(P),\quad\forall P\in\mathcal P
$$
\end{mydef}

\begin{mytheo}\rm  \label{teo4} Let $M_1:(\Omega,\mathcal A,\mathcal P)\pt\mathbb R^k$ and $M_2:(\Omega,\mathcal A,\mathcal P)\pt (\Omega_2,\mathcal A_2)$ be Markov kernels. If $M_2$ is sufficient, then there exists a regular conditional probability $P^{M_1|M_2}$ of $M_1$ given $M_2$ which is independent of $P\in\mathcal P$. There exists also a common version of the conditional expectations $E_P(M_1|M_2)$, $P\in\mathcal P$; it will be denoted $E(M_1|M_2)$. \end{mytheo}

 The next theorem extend to Markov kernels the Rao-Blackwell theorem.

\begin{mytheo}\rm \label{teo5} 
(Theorem of Rao-Blackwell generalized) Let $M_1:(\Omega,\mathcal
A,\mathcal P)\pt\mathbb R^k$ be an estimator of $f:\mathcal
P\rightarrow\mathbb R$ and $M_2:(\Omega,\mathcal A,\mathcal P)\pt
(\Omega_2,\mathcal A_2)$ be a sufficient Markov kernel for $\mathcal
P$. Then $E(M_1|M_2)$ is an estimator of $f$ with less convex risk
than $M_1$. If the loss function is strictly convex then, given
$P\in\mathcal P$, the risk at $P$ of $E(M_1|M_2)$ is strictly less
than the risk at $P$ of $M_1$ unless $E(M_1|M_2)\pi_2=\bar M_1$,
$P\otimes M_2$-a.s., where $\bar M_1$ is defined as in Theorem
\ref{teo2}. Finally, if $M_1$ is unbiased, so is $E(M_1|M_2)$.\end{mytheo}

\begin{myrem}\rm
Since $E(M_1|M_2)$ is a statistic, this theorem shows that the class
of non-randomized unbiased estimators of $f$ is complete in the
sense that, for every randomized unbiased estimator $M_1$ of $f$,
there exists a non-randomized unbiased estimator $E(M_1|M_2)$ with
less convex risk than $M_1$. Note that this assertion remains true
if the assumption of unbiasedness is dropped. This result
generalizes a similar result when $M_2$ is a  statistic rather than
a Markov kernel (for instance, see  Pfanzagl (1994, p. 105)).
\end{myrem}

\begin{mytheo}\rm \label{teo6}
(Theorem of Lehmann-Scheff\'{e} generalized) Let $M_1:(\Omega,\mathcal
A,\mathcal P)\pt\mathbb R^k$ be an unbiased estimator of $f:\mathcal
P\rightarrow\mathbb R^k$  and $M_2:(\Omega,\mathcal A,\mathcal P)\pt
(\Omega_2,\mathcal A_2)$ be a sufficient and complete Markov kernel
for $\mathcal P$. Then $E(M_1|M_2)$ is the estimator of $f$ which
minimizes the convex risk among all unbiased estimators of
$f$.\end{mytheo}

\section{Looking for an example of application of the generalized Lehmann-Scheffé Theorem}

As a Markov kernel is both an extension of the concepts of random variable and $\sigma$-field, Theorem \ref{teo6} can be considered as an unification of the statistics and $\sigma$-fields versions of the Lehmann-Sheffé Theorem.

But, beyond that, to give full meaning to this result, a sufficient and complete Markov kernel $M$ that is not associated to any statistic will be desirable. 
To this end, I wonder if the following general procedure can work: 
We start with a sufficient Markov kernel $M$  not associated to any statistic --two examples has been provided above-- and  construct a greater family of probabilities for which $M$ is complete and remains still sufficient. This way, if $M$ has finite mean $f$, the statistic $E(M|M)$  would be the estimator of $f$ which minimizes the convex risk among all unbiased estimators of $f$.

First, the right context is fixed. Let  $(\Omega,\mathcal A,\mathcal P)$ be a statistical experiment dominated by a $\sigma$-finite measure $\mu$. Let $P^*$ be a privileged dominating probability (i.e. $P^*$ is a probability measure on $(\Omega,\mathcal A)$ such that $\mathcal P\ll P^*$ and is of the form $P^*=\sum_n2^{-n}P_{\theta_n}$ for some countable dominating subfamily $\{P_{\theta_n}\colon n\ge 1\}$ of $\mathcal P$). A positive response to the following question would allow us to  resolve the posed problem, starting with the sufficient Markov kernel of Example \ref{continuousexample}.

\begin{myque}\rm \label{teo7} Let  $M:(\Omega,\mathcal A,\mathcal P)\pt (\Omega',\mathcal A')$ be a sufficient Markov kernel and denote $\hat{\mathcal P}$ the family of all probability measures $\hat P$ on $\mathcal A$ such that $\hat P\ll P^*$ and $M$ is sufficient for the extended family $\mathcal P\cup\{\hat P\}$. 
	
(a)  $M$ is sufficient and complete for $\hat{\mathcal P}$?

(b) Let $n,k\in\mathbb N$ and  suppose $(\Omega',\mathcal A')=(\mathbb R^n,\mathcal R^n)$. If $M$ has finite moments of order $k$ and $\hat{\mathcal P}_k$ denotes the subfamily of $\hat{\mathcal P}$ preserving this property of $M$, then $M$ is sufficient and complete for $\hat{\mathcal P}_k$?

	\end{myque}

Since, according to Heyer (1982, Theorem 22.3) (resp. Remark 1 of Definition \ref{def8}), a Markov kernel $M:(\Omega,\mathcal A,\mathcal P)\pt (\Omega',\mathcal A')$ is sufficient (resp. complete) if and only if the projection statistic 
$$\pi':(\omega,\omega')\in(\Omega\times\Omega',\mathcal A\times\mathcal A',\{P\otimes M\colon P\in\mathcal P\})\mapsto \omega'\in(\Omega',\mathcal A')$$
is sufficient (resp. complete), one might expect an answer to this question from the following --interesting in itself-- result, 


\begin{mytheo}\label{teon7} \rm  Let $S:(\Omega,\mathcal A,\mathcal P)\rightarrow(\Omega',\mathcal A')$ be a sufficient statistic and $\hat{\mathcal P}$ the family of all probabilities $\hat P$ on $(\Omega,\mathcal A)$  such that $\hat P\ll P^*$ and $S$ is sufficient for $\mathcal P\cup\{\hat P\}$. 

(a) Then $S$ is sufficient and complete for the  extended family $\hat{\mathcal P}$. 

(b) Suppose now that $(\Omega',\mathcal A')=(\mathbb R^k,\mathcal R^k)$, and $P^*\in\mathcal P$ (or there exists $P'\in\mathcal P$ such that $\mathcal P\ll P'$). If $S$ has finite moment of order $n\in\mathbb N$ and
$\hat{\mathcal P}_n$ denotes the set of the probability measures $\hat P\in\hat{\mathcal P}$ such that $S$ has finite $\hat P$-moment of order $n$. then $S$ is also sufficient and complete for $\hat{\mathcal P}_n$ .
	\end{mytheo}

Now we wonder if Question 1 actually becomes a consequence of the previous theorem. For that to be true, it would be enough that the extended family of $\{P\otimes M\colon P\in\mathcal P\}$ should coincide with $\{\hat P\otimes M\colon \hat P\in\hat{\mathcal P}\}$, 
something that does not seem very clear. 
So Question 1 remains unsolved but, in the attempt, Theorem \ref{teon7} has seen the light.

\section{Proofs}

{\sc Proof of Theorem \ref{teo1}.} 	First note that there exists a regular conditional probability $P^{M_1|M_2}$ (see Nogales (2013b)). It will be enough to show that, given $A_2\in\mathcal A_2$,
		$$\int_\Omega M_2(\omega,A_2)\int_{\mathbb R^k}xM_1(\omega,dx)dP(\omega)=\int_{A_2}  \int_{\mathbb R^k}xP^{M_1|M_2}(\omega_2,dx)dP^{M_2}(\omega_2)
		$$
		But by definition of $P^{M_1|M_2}$, for all $A_1,A_2$,
		$$\int_\Omega M_1(\omega,A_1) M_2(\omega,A_2)dP(\omega)=\int_{A_2}P^{M_1|M_2}(\omega_2,A_1)dP^{M_2}(\omega_2)$$
		i.e.,
		$$\int_\Omega M_2(\omega,A_2)\int_{\mathbb R^k}I_{A_1}(x)M_1(\omega,dx)dP(\omega)=\int_{A_2}  \int_{\mathbb R^k}I_{A_1}(x)P^{M_1|M_2}(\omega_2,dx)dP^{M_2}(\omega_2)$$
		It follows in a standard way that, for any nonnegative or $P^{M_1}$-integrable measurable function $f:\mathbb R^k\rightarrow\mathbb R^m$,
		$$\int_\Omega M_2(\omega,A_2)\int_{\mathbb R^k}f(x)M_1(\omega,dx)dP(\omega)=\int_{A_2}  \int_{\mathbb R^k}f(x)P^{M_1|M_2}(\omega_2,dx)dP^{M_2}(\omega_2)$$
		which gives the proof. $\Box$\vspace{2ex}

{\sc Proof of Theorem \ref{teo2}.}  Recall that $P^{M_2}=(P\otimes M_2)^{\pi_2}$.
			Now we define a Markov kernel $\hat
			M_1:(\Omega\times\Omega_2,\mathcal A\times\mathcal A_2)\pt\mathbb
			R^k$ by $\hat M_1((\omega,\omega_2),B)=M_1(\omega,B)$; $\hat M_1$ is
			the extension to $\Omega\times\Omega_2$ of $M_1$. We will prove that
			$(P\otimes M_2)^{\hat M_1|\pi_2}$ is a regular conditional
			$P$-probability of $M_1$ given $M_2$. We will use the following
			result from Nogales (2013b): ``If $T_2:(\Omega,\mathcal
			A)\rightarrow (\Omega_2,\mathcal A_2)$ is a random variable and
			$K_2(\omega,A_2)=\delta_{T_2(\omega)}(A_2)$ is its corresponding Markov
			kernel then, writing $P^{M_1|T_2}:=P^{M_1|K_2}$, we have
			$P^{M_1|T_2}(\cdot,A_1)=E_P(M_1(\cdot,A_1)|T_2).$''  Applying this
			result in the probability space $(\Omega\times\Omega_2,\mathcal
			A\times\mathcal A_2,P\otimes M_2)$, we have that, for
			$\omega_2\in\Omega_2$ and $B\in\mathcal R^k$,
			\begin{gather}(P\otimes M_2)^{\hat M_1|\pi_2}(\omega_2,B)=E_{P\otimes M_2}(\hat M_1(\cdot,B)|\pi_2=\omega_2)\label{2}
			\end{gather}
			
			Hence, given $A_2\in\mathcal A_2$,
			\begin{gather*}\begin{split}
					\int_{A_2}(P\otimes M_2)^{\hat M_1|\pi_2=\omega_2}(B)dP^{M_2}(\omega_2)&=\int_{A_2}(P\otimes M_2)^{\hat M_1|\pi_2=\omega_2}(B)d(P\otimes M_2)^{\pi_2}(\omega_2)\\&=
					\int_{A_2}E_{P\otimes M_2}(\hat M_1(\cdot,B)|\pi_2=\omega_2)d(P\otimes M_2)^{\pi_2}(\omega_2)\\&=
					\int_{\Omega\times A_2}M_1(\omega,B)d(P\otimes M_2)(\omega,\omega_2)\\&=
					\int_{\Omega}\int_{A_2}M_1(\omega,B)M_2(\omega,d\omega_2)dP(\omega)\\&=
					\int_{\Omega}M_1(\omega,B)M_2(\omega,A_2)dP(\omega)
				\end{split}\end{gather*}
				which proves that
				\begin{gather}
					(P\otimes M_2)^{\hat M_1|\pi_2}=P^{M_1|M_2}\label{1}
				\end{gather}
				
				Moreover, (\ref{2}) can be rewritten in the form
				$$\int_{\mathbb R^k}I_B(x)(P\otimes M_2)^{\hat M_1|\pi_2=\omega_2}(dx)=E_{P\otimes M_2}\left(\int_{\mathbb R^k}I_B(x)M_1(\cdot,dx)\left|\phantom{\int}\hspace{-1em}\right.\pi_2=\omega_2\right)
				$$
				It follows that, for a nonnegative or integrable measurable function $f:\mathbb R^k\rightarrow \mathbb R^m$,
				$$\int_{\mathbb R^k}f(x)(P\otimes M_2)^{\hat M_1|\pi_2=\omega_2}(dx)=E_{P\otimes M_2}\left(\int_{\mathbb R^k}f(x)M_1(\cdot,dx)\left|\phantom{\int}\hspace{-1em}\right.\pi_2=\omega_2\right)
				$$
				In particular, for $m=k$ and $f(x)=x$,
				$$\int_{\mathbb R^k}x(P\otimes M_2)^{\hat M_1|\pi_2=\omega_2}(dx)=E_{P\otimes M_2}\left(\int_{\mathbb R^k}x M_1(\cdot,dx)\left|\phantom{\int}\hspace{-1em}\right.\pi_2=\omega_2\right)
				$$
				Using (\ref{1}), we obtain
				$$E_P(M_1|M_2)=E_{P\otimes M_2}(\bar M_1|\pi_2).
				$$
				$\Box$\vspace{2ex}

{\sc Proof of Theorem \ref{teo3}.} 
{Let $T':(\Omega,\mathcal A)\rightarrow\mathbb R^k$ be} an arbitrary
unbiased estimator of $f$ and denote $\widetilde
T'(\omega,\omega_1):=T'(\omega)$. Hence $\widetilde T'$ is an unbiased
estimator of $f$ in the statistical experiment
$(\Omega\times\Omega_1, \mathcal A\times\mathcal A_1,\{P\otimes
M_1\colon P\in\mathcal P\})$. Since the coordinatewise projection
$\pi_1$ is sufficient, there exists a version of the conditional expectation $X'$ of $\widetilde T'$ given $\pi_1$ which is independent of $P\in\mathcal P$. The Rao-Blackwell theorem shows
that $X'\circ\pi_1$ has less covariance matrix than $\widetilde T'$.\par

Since $M_T=M_SM_1$, we have that, for all Borel set $B\in\mathcal R^k$ and all $\omega\in\Omega$,
$$I_B(T(\omega))=\int_{\Omega_1}I_B(S(\omega_1))M_1(\omega,d\omega_1)
$$
Hence, for all $\omega\in\Omega$, $S=T(\omega)$,
$M_1(\omega,\cdot)$-a.s. It follows that
$$\widetilde T(\omega,\omega_1)=(S\circ\pi_1)(\omega,\omega_1),\quad \{P\otimes M_1\colon P\in\mathcal P\}-\mbox{a.s.}
$$
where $\widetilde T(\omega,\omega_1)=T(\omega)$, for all $\omega\in\Omega$. So, $S$ is a conditional expectation of $\widetilde T$ given $\pi_1$ for all $P\in\mathcal P$.

The completeness of $\pi_1$ shows that $S\circ \pi_1=X'\circ\pi_1$,
$\{P\otimes M_1\colon P\in\mathcal P\}$-a.s., and this finish the
proof. $\Box$\vspace{2ex}

{\sc Proof of Theorem \ref{teo4}.} 
 According to Heyer (1982, Theorem 22.3), $M_2$ is sufficient {if and
 	only if} the coordinatewise projection $\pi_2:(\Omega\times\Omega_2,\mathcal A\times\mathcal A_2,\{P\otimes M_2\colon P\in\mathcal P\})\rightarrow (\Omega_2,\mathcal A_2)$ is sufficient. Landers and Rogge (1972, Theorem 7) shows the existence of a common regular conditional probability on $\mathcal R^k$ given $\pi_2$. The result follows from this fact and the following representation of the conditional distribution of $M_1$ given $M_2$ obtained in the proof of Theorem \ref{teo2}:
$$P^{M_1|M_2}(\omega_2,B)=(P\otimes M_2)^{\hat M_1|\pi_2}(\omega_2,B)=E_{P\otimes M_2}(\hat M_1(\cdot,B)|\pi_2=\omega_2)
$$
 The second assertion follows from this and Theorem \ref{teo1}.
$\Box$\vspace{2ex}

{\sc Proof of Theorem \ref{teo5}.}   $E(M_1|M_2)$ is well defined by the previous theorem and it is an unbiased estimator of $f$ by Corollary \ref{cor1}.
Moreover, if $W:\mathcal P\times\mathbb R^k\rightarrow[0,\infty[$ is
a convex loss function  (i.e., $W(P,\cdot)$ is a convex function for
every $P\in\mathcal P$) then applying the Jensen inequality (see
Pfanzagl (1994, Theorem 1.10.11)), we obtain from Theorem \ref{teo1}
that
\begin{gather*}\begin{split}
W(P,E_P(M_1|M_2))&=W\left(P,\int_{\mathbb R^k} xP^{M_1|M_2}(\cdot,dx)\right)\\
&\le\int_{\mathbb R^k}W(P,x)P^{M_1|M_2}(\cdot,dx)=E_P(W(P,M_1)|M_2),\quad P^{M_2}-\mbox{a.s.}
\end{split}\end{gather*}
where $W(P,M_1)$ denotes the kernel $W(P,\cdot)M_1$. The result
follows by integration with respect to $P^{M_2}$. Corollary
\ref{cor1} completes the proof in the unbiased case. $\Box$
\vspace{2ex}

{\sc Proof of Theorem \ref{teo6}.} 
If the Markov kernel $M'_1:(\Omega,\mathcal A,\mathcal P)\pt\mathbb
R^k$ is an arbitrary unbiased estimator of $f$  then, according to
the previous theorem, $X'_1:=E(M'_1|M_2)$ is a nonrandomized
unbiased estimator of $f$ with less convex risk than $M'_1$.
Moreover $X_1:=E(M_1|M_2)$ is an unbiased estimator of $f$; so
$E_{P^{M_2}}(X_1-X'_1)=0$ for all $P\in\mathcal P$. Since $M_2$ is
complete, we have that $X_1=X'_1$, $\{P^{M_2}\colon P\in\mathcal
P\}$-a.s. So $X_1$ has less convex risk than $M'_1$. 
$\Box$
\vspace{2ex}

\color{blue}
{\sc Proof of Theorem \ref{teon7}.} (a) By the sufficiency of $S$, $P^*$ being a privileged dominating probability, we can write the density of $P\in\mathcal P$ with respect to $P^*$ in the form $g_P\circ S$ for some suitable non-negative $\mathcal A'$-measurable function $g_P$. Such a factorization is also valid for $\hat P\in\hat{\mathcal P}$, as these properties of $S$ and $P^*$ remain valid for $\mathcal P\cup\{\hat P\}$. Therefore, $S$ is sufficient for the extended family $\hat{\mathcal P}$.
	
	To prove the $\hat{\mathcal P}$-completeness of $S$, a reasoning by reduction to the absurd will be made by assuming the existence of a probability $\hat P_0\in\hat{\mathcal P}$ and a random variable $f:(\Omega',\mathcal A')\rightarrow\mathbb R$ such that
	$$E_{\hat P^S}(f)=0,\ \ \forall \hat P\in\hat{\mathcal P},\qquad\text{and}\qquad
	\hat P_0^S(f\ne 0)>0
	$$
	Without loss of generality we can suppose that $\hat P_0^S(f>0)>0$. Since $\hat P_0^S\ll {P^*}^S$, $\alpha:={P^*}^S(f>0)>0$. So $g_1:=\frac1{\alpha}I_{\{f>0\}}$ defines a ${P^*}^S$-density of a probability measure on $\mathcal A'$ and $g_1\circ S$ becomes a $P^*$-density of a probability measure $P_1$ on $\mathcal A$. So $P_1\in\hat{\mathcal P}$ and $dP_1^S=g_1d{P^*}^S$. Moreover,
	$$E_{P_1^S}(f)=\frac1{\alpha}\int_{\Omega'}f\cdot I_{\{f>0\}}d{P^*}^S>0,
	$$
	a contradiction which finishes the proof of (a).
	
	(b) The sufficiency is obtained in the same way as in (a). To prove completeness, taking $k=1$ for simplicity, we just have to check that $P_1\in\hat{\mathcal P}_n$. But
	$$E_{P_1}(|S|^n)=\frac1{\alpha}\int_{\Omega}|S(\omega)|^n\cdot I_{\{f>0\}}(S(\omega))dP^*(\omega)\le\frac1{\alpha}E_{P^*}(|S|^n)<\infty.
	\quad\Box$$
\vspace{2ex}

\color{black}

\section* {References:}

\begin{itemize}
	\item 
Blackwell, D.:  Conditional expectation and unbiased sequential estimation, Ann. Math. Statist. 18, 105-110 (1947).

	\item Dellacherie, C., Meyer, P.A.:  Probabilities and Potentiel C, North-Holland, Amsterdam (1988).

	\item Heyer, H.:  Theory of Statistical Experiments, Springer, Berlin (1982).

    \item Kolmogorov, A.N.:  Grundbegriffe der Wahrscheinlichkeitsrechnung, Sprin\-ger, Berlin (1933).

	\item Landers, D. and  Rogge, L.: A note on completeness, Scandinavian J. Statist. 3, 139 (1976).

	\item Lehmann, E.L., Scheff\'{e}, H.: Completeness, similar regions, and unbiased estimation, Sankhy$\bar{\mbox{a}}$ 10, 305-340; 15, 219-236; Correction 17, 250 (1950, 1955, 1956).

	\item Nogales, A.G.: On Independence of Markov Kernels and a Generalization of Two Theorems of Basu, Journal of Statistical Planning and Inference 143, 603-610  (2013a).

	\item Nogales, A.G.: Existence of Regular Conditional Probabilities for Markov kernels, Statistics and Probability Letters 83, 891-897 (2013b).

	\item Pfanzagl, J.: Parametric Statistical Theory, de Gruyter, Berlin (1994).

	\item Rao, C.R.:  Information and the accuracy attainable in the estimation of statistical experiments, Bull. Calcutta Math. Soc. 37, 81-91 (1945).

\end{itemize}

\end{document}